# Icosahedral Tiling with Dodecahedral Structures

Mehmet Koca[a], Ramazan Koc[b], Nazife Ozdes Koca[c,*] and Abeer Al-Siyabi[c]

[a]Department of Physics, Cukurova University, Adana, Turkey, retired
[b]Department of Physics, Gaziantep University, Gaziantep, Turkey
[c]Department of Physics, College of Science, Sultan Qaboos University, P.O. Box 36, Al-Khoud, 123 Muscat, Sultanate of Oman, *Correspondence e-mail: nazife@squ.edu.om

## ABSTRACT

Icosahedron and dodecahedron can be dissected into tetrahedral tiles projected from 3D-facets of the Delone polytopes representing the deep and shallow holes of the root lattice $D_6$. The six fundamental tiles of tetrahedra of edge lengths 1 and $\tau$ are assembled into four composite tiles whose faces are normal to the 5-fold axes of the icosahedral group. The 3D Euclidean space is tiled face-to-face by the composite tiles with an inflation factor $\tau$ generated by an inflation matrix. The aperiodic tiling is a generalization of the Tubingen triangular tiling in 2-dimensions for the faces of the tiles are made of Robinson triangles. Certain combinations of the tiles constitute dodecahedra with edge lengths of 1 and the golden ratio $\tau = \frac{1+\sqrt{5}}{2}$.

# 1. Introduction

Icosahedral quasicrystallography is the focal interests of many scientists from diverse fields of research. For a review see for instance the references (Di Vincenzo & Steinhardt, 1991; Janot, 1993; Senechal, 1995). The subject is mathematically challenging as it requires the aperiodic tiling of the space by some prototiles. For a an excellent exposition we propose the references (Baake & Grimm, 2013; Baake & Grimm, 2020).

There have been two intimately related major approaches for the aperiodic order of the 3D space with icosahedral symmetry. In the increasing order of symmetry the set of four tiles (Socolar & Steinhardt, 1986) consists of acute rhombohedron, Bilinski rhombic dodecahedron, rhombic icosahedron and rhombic triacontahedron; the latter three are constructed with two Ammann tiles of acute and obtuse rhombohedra. Decorations of the Ammann tiles were proposed by Katz (Katz, 1989) and recently reviewed by Hann-Socolar-Steinhardt (Hann, Socolar & Steinhardt, 2018). More fundamental tiles were proposed by Danzer known as *ABCK* tetrahedral tiling (Danzer, 1989). Later, it was shown that these two sets of tiles are related to each other (Danzer, Papadopolos & Talis, 1993; Roth, 1993) as both sets of tiles are composed of faces normal to the 2-fold axes of the icosahedral symmetry. Ammann rhombohedral and Danzer *ABCK* tetrahedral tilings can be obtained from the projections of the six dimensional cubic lattice $B_6$ and the root lattice $D_6$ respectively (Koca, Koca & Koc, 2015; Al-Siyabi, Koca & Koca, 2020). Kramer and Andrle (Kramer & Andrle, 2004) have also investigated the Danzer tiles in the context of $D_6$ lattice with its relation to the wavelets.

For the icosahedral aperiodic order it is certainly desirable to search for the prototiles with faces normal to the 5-fold and/or 3-fold axes. It is well known that (Conway & Sloane, 1999) the Delone polytopes of the the root lattice $D_6$ defined by the weight vector $\omega_1$ representing the cross polytope with 12 vertices and the weight vectors $\omega_5$ and $\omega_6$ each representing a hemi-cube with 32 vertices (Coxeter, 1973) tile the the root lattice $D_6$ in an alternating order by centralizing the vertices of the Voronoi cell of the root lattice. Projection of the Delone cell represented by $\omega_1$ form an icosahedron dissected into four types of tetrahedral tiles obtained from the projections of the 240 3D-facets of the cross polytope. Similarly the hemi-cubes project into icosahedra and dodecahedra where the dodecahedra are dissected into six tetrahedral tiles (including former four tiles) projected from the 640+640 3D-facets of the hemi-cubes. We call them the *fundamental* tiles faces of which are composed of triangles normal to the 5-fold and 3-fold axes. As the faces of the fundamental tiles being normal to the 3-fold axes they cannot be partitioned in terms of similar equilateral triangles. We define a new set of four *composite* prototiles assembled by the fundamental tiles whose faces are normal to the 5-fold axes only as such they are made of Robinson triangles. The composite tiles can then be inflated by an inflation factor $\tau$ which can be generated by an inflation matrix. In what follows we describe the procedure how to tile the 3D space with four composite tiles.

The paper is organized as follows. In Sec. 2, we introduce the fundamental tiles projected from the Delone cells of the root lattice $D_6$ which admits the icosahedral group $H_3$ as a maximal subgroup. Compositions of the icosahedra and dodecahedra in terms of fundamental tiles will be briefly discussed. In Sec 3 we introduce four composite tiles assembled by the fundamental tiles so that they are composed of the Robinson triangles normal to the 5-fold axes. The composite tiles can be inflated with the $4 \times 4$ inflation matrix whose eigenvalues are $\tau^3, \tau, \sigma$ and $\sigma^3$ where $\sigma = -\tau^{-1} = \frac{1-\sqrt{5}}{2}$ is the algebraic conjugate of $\tau$. Right and left eigenvectors of the Perron-Frobenius (PF) eigenvalue $\tau^3$ and the projection matrix are calculated. The composite tiles involve dodecahedral structures of edge lengths 1 and $\tau$ already in the 2$^{nd}$ and 3$^{rd}$ order of the inflation. As the inflation gets larger and larger the dodecahedra pup up at certain points of the space and the numbers of dodecahedra and tiles connecting them can be determined from the inflation matrix. For the faces of the composite tiles are made of Robinson triangles its relevance to the Penrose-Robinson tiling (PRT) and possible relations to the Tubingen triangle tiling (TTT) is briefly mentioned in the concluding remarks.

## 2. Projections of the 3D facets of Delone cells of the root lattice $D_6$

The Delone cell characterized by the weight vector $\omega_1$ possesses 12 vertices, 60 edges, 120 triangular faces, 240 tetrahedral facets, 192 4-simplexes and 64 5-simplexes. The orthogonal projection of the polytope into 3D space defined as $E_\parallel$ form an icosahedron whose constituents are four types of tetrahedra $t_1, t_2, t_5$ and $t_6$. Similarly each of the hemi-cubes has 32 vertices, 240 edges, 640 triangular faces, 640 tetrahedral facets, 252 4-simpexes and 44 5-simplexes. They project into icosahedron or dodecahedron with six different tiles also including tiles $t_3$ and $t_4$ as shown in Table 1.

**Table 1**
The fundamental tiles projected from Delone cells of $D_6$.

| Name of tile | Sketch | Number of faces $(a, b, c)$ | Volume |
|---|---|---|---|
| $t_1$ | 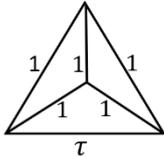 | $2 \times (1,1,1)$<br>$2 \times (1,1,\tau)$ | $\dfrac{1}{12}$ |
| $t_2$ | 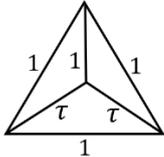 | $1 \times (1,1,1)$<br>$2 \times (1,1,\tau)$<br>$1 \times (1,\tau,\tau)$ | $\dfrac{\tau}{12}$ |
| $t_3$ | 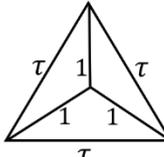 | $1 \times (\tau,\tau,\tau)$<br>$3 \times (1,1,\tau)$ | $\dfrac{\tau}{12}$ |
| $t_4$ | 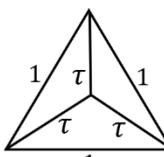 | $1 \times (1,1,1)$<br>$3 \times (1,\tau,\tau)$ | $\dfrac{\tau^2}{12}$ |
| $t_5$ | 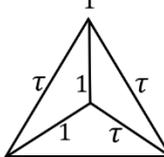 | $1 \times (\tau,\tau,\tau)$<br>$1 \times (1,1,\tau)$<br>$2 \times (1,\tau,\tau)$ | $\dfrac{\tau^2}{12}$ |
| $t_6$ | 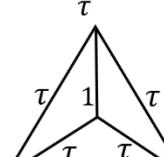 | $2 \times (\tau,\tau,\tau)$<br>$2 \times (1,\tau,\tau)$ | $\dfrac{\tau^3}{12}$ |

The edge lengths of six tetrahedral tiles are either 1 or $\tau$ besides an overall factor arising from the projection. Tile $t_1$ consists of two equilateral triangles of edge lengths 1 and two Robinson triangles of edge lengths $(1,1,\tau)$. The tile $t_2$ is made of one equilateral triangle of edge length 1, two Robinson triangles of edge lengths $(1,1,\tau)$ and one Robinson triangle with edge lengths $(1,\tau,\tau)$. Tile $t_3$ has one triangular face of edge lengths $\tau$ and the others are Robinson triangles of $(1,1,\tau)$. The set of faces of tile $t_4$ consists of one triangle with edge lengths 1 and the others are the Robinson triangles with edges $(1,\tau,\tau)$. The tile $t_5$ consists of one equilateral triangle with edge length $\tau$, two Robinson triangles with edges $(1,\tau,\tau)$ and one triangle with edges $(1,1,\tau)$. Lastly, tile $t_6$ has 2 faces with edges $(1,\tau,\tau)$ and two equilateral faces with edge lengths $\tau$. All Robinson triangles and equilateral triangles are normal to 5-fold and 3-fold axes respectively.

Let an icosahedron of edge length $\tau^n$ is denoted by $i(\tau^n), (n = 0, 1)$ then their constituents in terms of fundamental tiles are given by

$$i(1) = 7t_1 + 6t_2 + 2t_5 + t_6,$$
$$i(\tau) = t_1 + 8t_2 + 10t_3 + 10t_4 + 16t_5 + 3t_6. \qquad (1)$$

There are two ways of constructions of icosahedron $i(1)$: either constructing via a pentagonal antiprism and then adding two pentagonal pyramids or first forming a three diminished icosahedron (Johnson solid $J_{63}$) and then adding three pentagonal pyramids. We will explain the construction of the second kind following even a simpler technique. First, assemble three tiles $2t_5 + t_6$ into a composite tile (denoted by $T_3 =: t_5 + t_6 + t_5$ for further use) by matching their equilateral triangular faces placing $t_6$ in between two tiles of $t_5$. They form a "pentagonal pyramid" with pentagonal base of edge length 1 and five Robinson triangles with edges $(1,\tau,\tau)$. This polyhedron was used in a 7-tile system by Kramer (Kramer, 1982) and it was called a tent but we will continue naming it as "pentagonal pyramid" although it is not actual pentagonal pyramid. Then we glue 5 tiles of $t_2$ matching the Robinson triangles of $T_3$, then fill the gaps between $t_2$ with 5 tiles of $t_1$. It is almost done except covering the pentagonal base of the "pentagonal pyramid" by a pentagonal pyramid formed by a sandwich of tiles $t_1 + t_2 + t_1$. The result is an icosahedron $i(1)$. What is really happenning is that the cross polytope actually collapses into an icosahedron with all its facets converting themselves to four tiles described above. We leave the construction of the icosahedron $i(\tau)$ to the reader.

Let the dodecahedron of edge length $\tau^n$ will be denoted by $d(\tau^n)$. We now describe the construction of dodecahedron $d(1)$ in terms of six fundamental tiles. The tile content of dodecahedron $d(1)$ is given by

$$d(1) = 3t_1 + 4t_2 + 10t_3 + 10t_4 + 4t_5 + 7t_6. \qquad (2)$$

This construction allows us to define the composite tiles, therefore, we will defer its construction to the next section.

## 3. Composite tiles and the inflation matrix

The composition of the fundamental tiles in dodecahedron $d(1)$ in (2) is such that they allow us to define the following composite tiles. Since dodecahedron has faces normal to the 5-fold axes then the tiles constituting it must have faces consisting of Robinson triangles only. This allows us to compose the new tiles by matching the equilateral triangular faces of the fundamental tiles. What we see in the construction of $d(1)$ is the formation of composite tiles on their equilateral triangular faces as follows:

$$T_1 =: E + C; \ E =: t_4 + t_1 + t_4, \ C =: t_3 + t_6 + t_3,$$

$$T_2 =: t_2 + t_4, \qquad (3)$$

$$T_3 =: t_5 + t_6 + t_5,$$

$$T_4 =: t_3 + t_6 + t_5.$$

They are constructed in such a way that the the composite tiles have no reference to the equilateral triangular faces of the fundamental tiles but they dispay the faces composed of Robinson triangles. The composite tiles are illustrated in Table 2 with their geometric structures.

**Table 2**
The composite tiles ($N_0$: number of vertices, $N_1$: number of edges, $N_2$: number of faces)

| Name of tile | Figure | $N_0$ | $N_1$ | $N_2$ | Type of faces | Volume |
|---|---|---|---|---|---|---|
| $T_1$ | 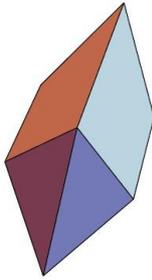 | 8 | 14 | 8 | $4 \times (1,1,\tau)$<br>$4 \times (1,1,1,\tau)$ (trapezoid) | $\dfrac{2\tau^4}{12}$ |
| $T_2$ | 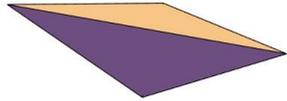 | 4 | 6 | 4 | $2 \times (1,\tau,\tau)$<br>$2 \times (\tau^2,\tau,\tau)$ | $\dfrac{\tau^3}{12}$ |
| $T_3$ | 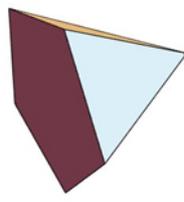 | 6 | 10 | 6 | $5 \times (1,\tau,\tau)$<br>1 pentagon of edge length 1 | $\dfrac{4\tau+3}{12}$ |
| $T_4$ | 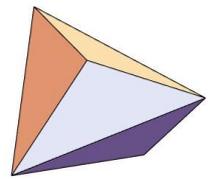 | 6 | 11 | 7 | $3 \times (1,1,\tau)$<br>$3 \times (1,\tau,\tau)$<br>$1(1,1,1,\tau)$ (trapezoid) | $\dfrac{2\tau^3}{12}$ |

The tile $T_1 = E + C$ is made of two composite tiles as it always occurs in this combinations in any $d(\tau^n)$. The tile $E$ is a nonconvex octahedron obtained by matching equilateral triangular faces of two $t_4$ with the two equilateral triangular faces of $t_1$. It has $N_0 = 6$ vertices, $N_1 = 12$ edges and $N_2 = 8$ faces. The set of faces of the composite tile $E$ is $6 \times (1, \tau, \tau)$ and $2 \times (1,1, \tau)$. The tile $C$ is composed by sandwiching a $t_6$ between two $t_3$ tiles on their equilateral triangular faces and has $N_0 = 6$ vertices, $N_1 = 12$ edges and $N_2 = 8$ faces same as $E$ but the number of faces are exchanged as $6 \times (1,1, \tau)$ and $2 \times (1, \tau, \tau)$. The tile $T_1$ is obtained by inserting $C$ between the legs of $E$ by matching two faces $(1, \tau, \tau)$ with similar faces of $E$. The new composite tile $T_1$ consists of $N_0 = 8$ vertices, $N_1 = 14$ edges and $N_2 = 8$ faces. The tile $T_1$ has $4 \times (1,1, \tau)$ triangles and $4 \times (1,1,1, \tau)$ quadrilaterals. The last four faces are isosceles trapezoids made of $(1,1, \tau)$ and $(1, \tau, \tau)$ triangles.

The tile $T_2$ which was also used in a 7-tile system of Kramer (Kramer, 1982) is a tetrahedron with faces $2 \times (1, \tau, \tau)$ and $2 \times (\tau^2, \tau, \tau)$ which is obtained by gluing two equilateral faces of tiles $t_2$ and $t_4$.

The tile $T_3$ is already described in Sec.2 which has 6 vertices, 10 edges and 6 faces, one is a pentagon the other 5 are of type $(1, \tau, \tau)$. There is another version of $T_3$ which can be denoted by $t_5 + t_6 + t^5$, where one of the $t_5$ is rotated so that one obtains a face of type $(1,1,1, \tau)$ instead of a pentagon. This composite tile is required in the construction of icosahedron if one starts with the pentagonal antiprism.

The tile $T_4$ is obtained from $T_3$ by replacing one of $t_5$ by $t_3$. Further properties of the composite tiles can be obtained from Table 2. Just to mention another common propery is that the dihedral angles between faces of the composite tiles are either $\tan^{-1}(2)$ or $\pi - \tan^{-1}(2)$. It is a matter of exercise how to build dodecahedron $d(1)$ in terms of composite tiles. It consists of two parts, one being twice the other volume-wise; the smaller part consists of $T_1 + 2T_2 + T_4$ assembled face-to-face matching and the larger part is built as $2T_1 + 2T_2 + 3T_4$. These two polyhedra form frustrums and can be matched at their pentagonal faces of edge length $\tau$ leading to the dodecahedron

$$d(1) = 3T_1 + 4T_2 + 4T_4. \tag{4}$$

The tiles are combined in such a way that they meet at three points inside dodecahedron $d(1)$ forming the vertices of an equilteral triangle of edge length 1. The dodecahedron $d(\tau)$ can be constructed from the icosahedron $i(1)$ given in (1) by covering the equilateral faces of $t_1$ and $t_2$ with the equilateral faces of $t_4$. By this one notes that the polyhedron is a star icosahedron and with this construction the tiles $t_1$ and $t_2$ are converted to the composite tiles $E$ and $T_2$. Filling the gaps between the legs of $E$ by the tiles $C$ one obtains $7T_1$ composite tiles. The rest follows by face-to-face matching to complete the construction of the dodecahedron given by

$$d(\tau) = 7T_1 + 18T_2 + 14T_3 + 10T_4. \tag{5}$$

Before we proceed further we should mention that there are 12 vertices of icosahedron in the dodecahedron $d(\tau)$ with no face structures as they are covered by the tiles $t_4$. One can infer an inflation rule with an inflation factor $\tau$ by comparing (4) and (5):

$$d(\tau) = 3\tau T_1 + 4\tau T_2 + 4\tau T_4,$$

where $\tau T_1 = \tau E + \tau C$ with $\tau E = 2T_2 + T_3 + T_4$ and $\tau C = T_1 + T_3 + T_4$ so that

$$\tau T_1 = T_1 + 2T_2 + 2T_3 + 2T_4. \tag{6}$$

Inflation of the other composite tiles can be constructed as follows

$$\tau T_2 = 2T_2 + T_3,$$
$$\tau T_3 = T_1 + 2T_2 + T_3 + T_4,$$
$$\tau T_4 = T_1 + T_2 + T_3 + T_4, \qquad (7)$$

where $\tau T_i$ $(i = 1, 2, 3, 4)$ are illustrated in Fig.1.

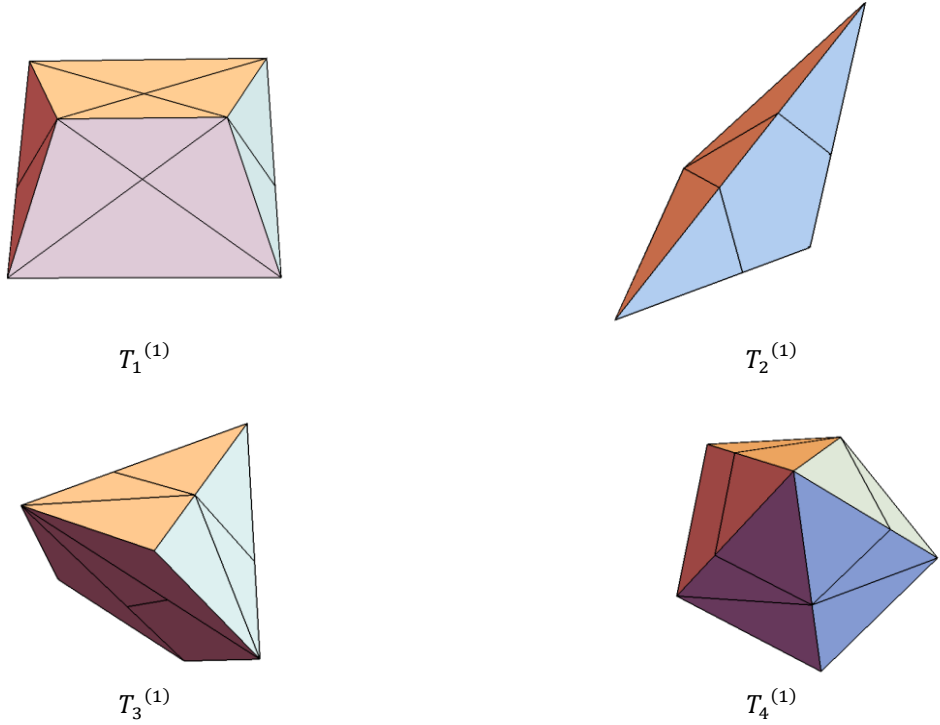

**Figure 1**
Composite tiles obtained from fundamental tiles inflated by the factor $\tau$ (see equation 10 for the definition).

The relations in (6-7) can be combined in a matrix equation,

$$\tau T_i = \sum_{j=1}^{4} M_{ij}\, T_j, \ (i = 1, 2, 3, 4)$$

where the matrix $M$ can be written as

$$M = \begin{pmatrix} 1 & 2 & 2 & 2 \\ 0 & 2 & 1 & 0 \\ 1 & 2 & 1 & 1 \\ 1 & 1 & 1 & 1 \end{pmatrix}. \qquad (8)$$

The eigenvalues of the inflation matrix are $\tau^3$, $\tau$, $\sigma$ and $\sigma^3$; the right eigenvector corresponding to the Perron-Frobenius (PF) eigenvalue $\tau^3$ has the components $(V_{T_1}, V_{T_2}, V_{T_3}, V_{T_4})^T$ with statistical normalization it reads $(\sigma^2, -\frac{\sigma^3}{2}, \frac{4\sigma+3}{2}, -\sigma^3)^T \cong (0.3820, 0.1180, 0.2639, 0.2361)^T$. This implies that the tile $T_1$ occupies most volume of the aperiodic tiling, nearly 38% of the space. The statistically normalized left eigenvector or the right eigenvector of $M^T$ of the inflation matrix is $\frac{2}{5\tau+4}(\frac{\tau}{2}, \tau^2, \tau, 1)^T \cong (0.1338, 0.4331, 0.2677, 0.1654)^T$ and it shows the relative frequency of the tiles indicating that the tile $T_2$ is nearly 43% more frequent. The PF projection matrix is determined as

$$\lim_{n\to\infty} \tau^{-3n} M^n = P = \frac{1}{30}\begin{pmatrix} 2(\tau+2) & 4(3\tau+1) & 4(\tau+2) & 4\sqrt{5} \\ \sqrt{5} & 2(\tau+2) & 2\sqrt{5} & 2(2+\sigma) \\ 5 & 10\tau & 10 & -10\sigma \\ 2\sqrt{5} & 4(\tau+2) & 4\sqrt{5} & 4(2+\sigma) \end{pmatrix}, \quad P^2 = P. \quad (9)$$

After this general procedure we will illustrate some of the inflated patches. For this, we first define the inflated tiles by a new notation. Let us denote by $\tau^n T_i =: T_i^{(n)}, n = 0,1,2,...$ then we can write

$$T_i^{(n)} = \sum_{j=1}^{4}(M^n)_{ij} T_j, \quad (i = 1, 2, 3, 4). \quad (10)$$

Certain dodecahedra of types $d(1)$ and $d(\tau)$ can be constructed when the higher order inflation factor is applied on the composite tiles. For example, we obtain the dodecahedron $d(1)$ in the inflated tiles given by

$$\begin{aligned} T_1^{(2)} &= d(1) + 2T_2^{(1)} + T_3^{(1)} + T_4^{(1)} + T_2 + 4T_3, \\ T_2^{(3)} &= d(1) + 2T_2^{(2)} + 5T_2 + 6T_3, \\ T_3^{(2)} &= d(1) + 5T_2 + 6T_3, \\ T_4^{(2)} &= d(1) + 3T_2 + 5T_3, \end{aligned} \quad (11)$$

where they are depicted in Fig. 2.

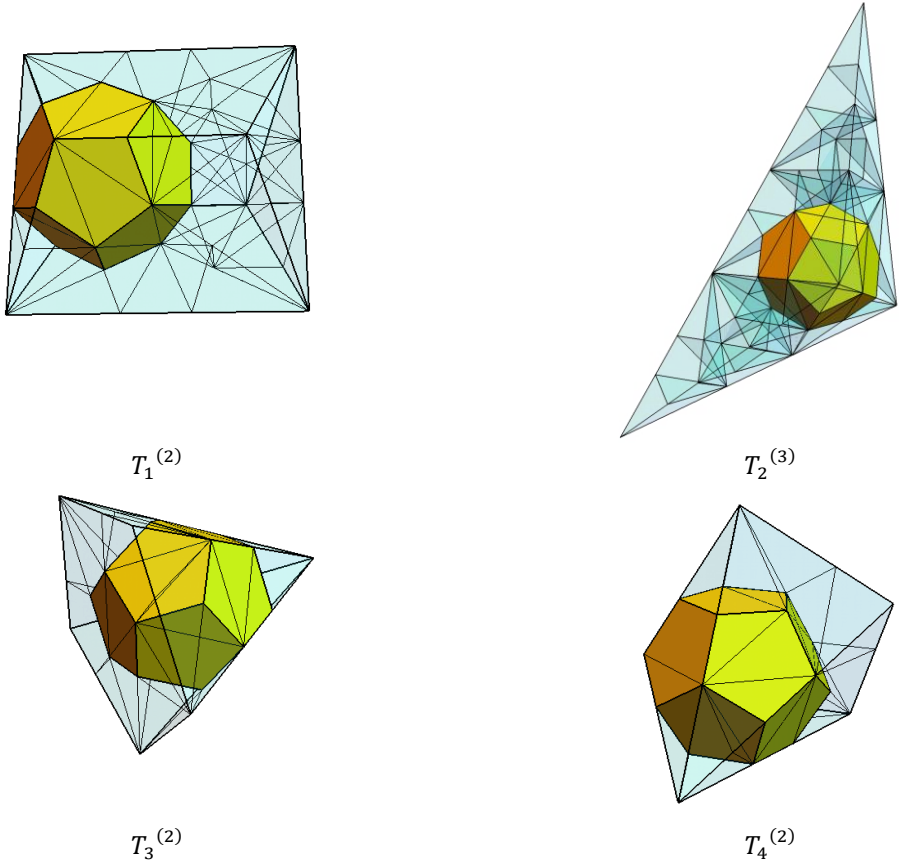

**Figure 2**
Illustration of tiles in (11).

Further inflation of tiles in (11) by $\tau$ will produce $d(\tau)$ in three of the inflated tiles while we obtain in

$$T_2^{(4)} = 2d(1) + d(\tau) + 4T_2^{(2)} + 5T_2^{(1)} + 6T_3^{(1)} + 10T_2 + 12T_3, \qquad (12)$$

where both dodecahedra $d(1)$ and $d(\tau)$ occur similtaneously. Another interesting case happens in the inflation represented by

$$T_1^{(4)} = 13d(1) + 2d(\tau) + 9T_2^{(2)} + 14T_2^{(1)} + 14\,T_3^{(1)} + 3T_4^{(1)} + 45\,T_2 + 68T_3. \qquad (13)$$

Similar formulae can be obtained for $T_3^{(4)}$ and $T_4^{(4)}$. They are illustrated in Fig. 3 by highlighting the dodecahedral structures and leaving the others transparent.

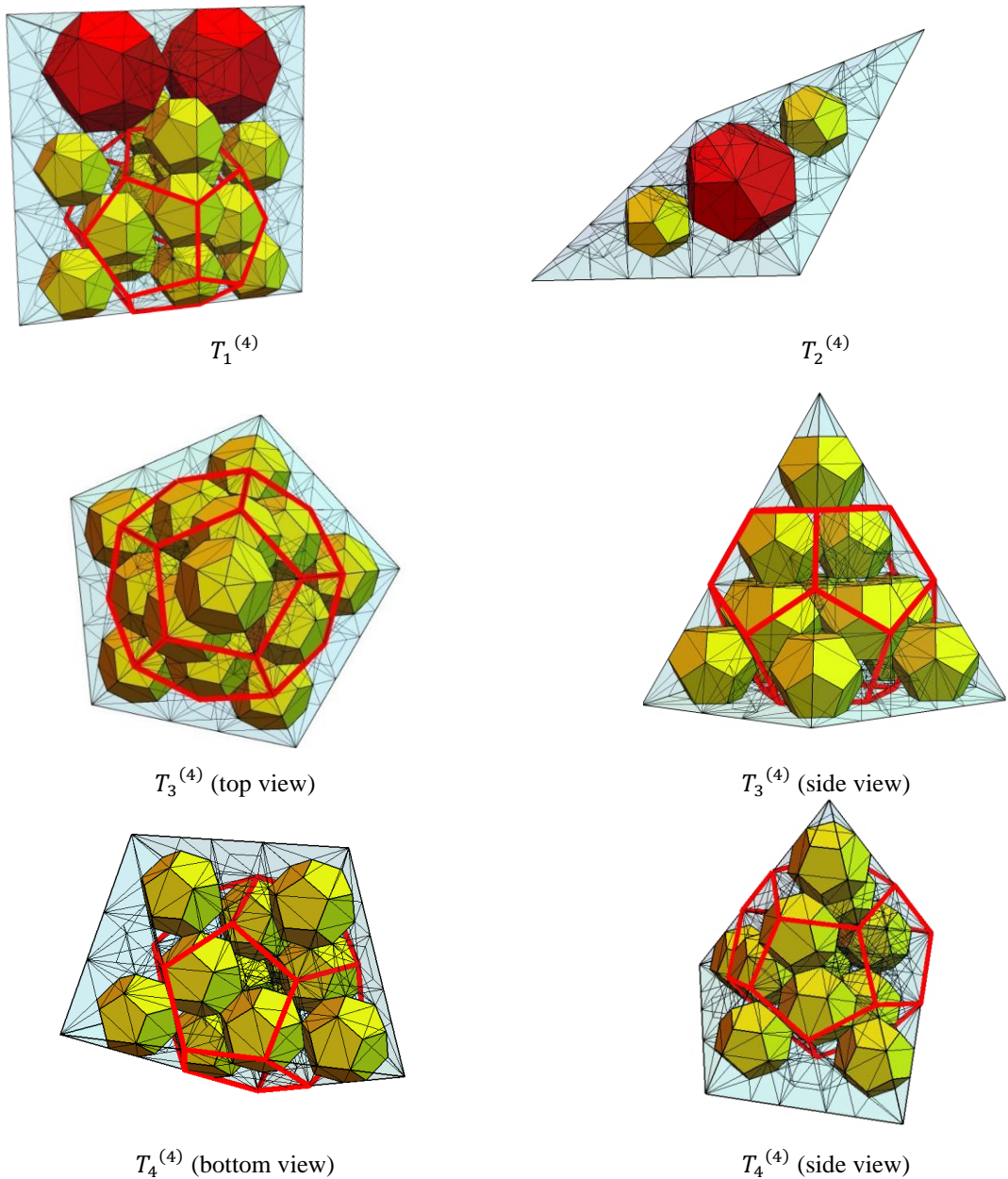

$T_1^{(4)}$         $T_2^{(4)}$

$T_3^{(4)}$ (top view)         $T_3^{(4)}$ (side view)

$T_4^{(4)}$ (bottom view)         $T_4^{(4)}$ (side view)

**Figure 3**
An illustration of $T_i^{(4)} (i = 1, 2, 3, 4)$ where the tiles $d(1)$ and $d(\tau)$ are demonstrated in different colours and dodecahedral frames indicate $d(\tau^2)$.

The dodecahedron $d(1)$ can be inflated to an arbitrary order of the inflation factor $\tau$. They all reduce to a number of $d(1)$ and/or $d(\tau)$ along with the other composite tiles. In Fig. 4 we illustrate $d(\tau^2)$ which is composed of $7d(1)$ and the accompanying composite tiles.

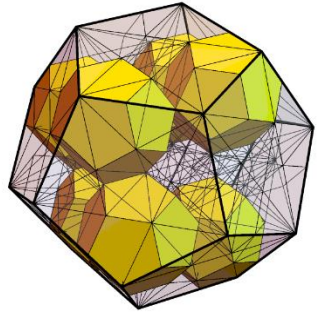

**Figure 4**
Dissection of $d(\tau^2)$ (composite tiles are transparent).

It is clear that the patches include $d(1)$ and $d(\tau)$ in abundance. Every dodecahedron $d(\tau^n), n \geq 2$ can be dissected into dodecahedra $d(1)$ and $d(\tau)$ along with other composite tiles. To give an example consider the dodecahedron $d(1)$ inflated by $\tau^{10}$,

$$d(\tau^{10}) = 432139 d(1) + 92850 d(\tau) + 1064050 T_1 + 6341550 T_2 + 4720730 T_3 + 1064050 T_4. \quad (14)$$

It is clear how fast the number of dodecahedra $d(1)$ and $d(\tau)$ are growing which are connected by the composite tiles.

### 4. Concluding remarks

The present tiling scheme is an alternative model to the widely known *ABCK* Danzer tiling which also involves the underlying principles of the Ammann tiling with acute and obtuse rhombohedra whose faces are normal to the 2-fold axes of the icosahedral symmetry. Our model differs from the *ABCK* tiling because the faces of the tiles are normal to the 5-fold axes and, not only this, the faces of the tiles are composed of the Robinson triangles in a manner of the Tubingen triangle tiling (TTT) which possesses 5-fold planar symmetry provided the original faces of the fundamental tile are preserved. In this sense the model is a generalization of the TTT to 3D space with icosahedral symmetry with dodecahedral structures. Modification of the model seems to be possible as a generatization of the Penrose-Robinson tiling (PRT) or mixture of two tilings by replacing the tile $T_3$ by $\overline{T_3} = t_5 + t_6 + t^5$. A detailed exposure of the work including its relevance to the lattice $D_6$ will be submitted for publication.